\documentclass[12pt]{amsart}
\usepackage{amsmath,amssymb,amscd,amsthm}
\usepackage[all,cmtip]{xy} 
\usepackage{cite}

\usepackage[a4paper,divide={2cm,*,2cm},nofoot]{geometry} 
\usepackage[all]{xypic}
\usepackage{graphicx}
\usepackage{mathrsfs}

\newtheorem{thm}{Theorem}[section]

\newtheorem{lem}[thm]{Lemma}

\theoremstyle{definition}
\newtheorem{defn}[thm]{Definition}
\theoremstyle{remark}
\newtheorem{rem}[thm]{Remark}
\numberwithin{equation}{section}

\newcommand{\C}{\mathbb{C}} 
\newcommand{\R}{\mathbb{R}}

\DeclareMathOperator{\Sym}{Sym}
\DeclareMathOperator{\Herm}{Herm}

\DeclareMathOperator{\tr}{tr}

\providecommand{\norm}[1]{\lVert#1\rVert}

\usepackage{esint}

\newcommand{\El}{\mathcal{E}}
\newcommand{\wt}[1]{\widetilde{#1}}
\newcommand{\proj}{\mathrm{p}}

\title{A Remark On $C^{2,\alpha}$-Regularity of the Complex Monge-Amp\`ere Equation}
\author{Yu Wang}

\begin{document}

\begin{abstract}
We prove the $C^{2,\alpha}$-regularity of the solution $u$ of the equation
\[
\det( u_{\bar{k} j}) = f,  \quad  f^{1/n} \in  C^{\alpha}, \quad  f \geq \lambda
\]
under the assumption in upper bound of $\Delta u$. Our result settles down the regularity problem related to the paper \cite{Tian} (also see \cite{Zhang}). 
\end{abstract}

\maketitle


\section{Introduction}
In this note, we prove 
\begin{thm}
\label{CMA}
Let $u \in C^2 (B_1)$ be a plurisubharmonic function that solves the equation
\begin{equation}
\label{CMA}
\det (u_{\bar{k} j}) = f \; \text{ in } B_1.
\end{equation}
Suppose $0 < \alpha < 1$
\begin{equation}
\label{c1}
f^{1/n} \in C^{\alpha} (\overline{B_1})  ,\quad   \inf_{B_1} f^{1/n} \geq \lambda > 0
\end{equation}
and 
\begin{equation}
\label{c2}
 \sup_{B_1 } \Delta u  \leq \Lambda.
\end{equation}
Then, for any $0 < \beta < \alpha$, there exists a constant $C$ depending on $\beta, n, \lambda, \Lambda, \norm{f^{1//n}}_{C^{\alpha} (\overline{B_1})}, \norm{u}_{L^{\infty}(B_1)}$ such that
\[
u \in C^{2, \beta} (\overline{B_{1/2}}) \; \text{ and } \;  \norm{u}_{C^{2,\beta}(\overline{B_{1/2}})} \leq C.
\]
\end{thm}

By the standard nonlinear elliptic theory (Thm.3 in \cite{Caf} and Thm.6.6 in \cite{CC}), the above regularity result is a direct consequence of the following theorem:

\begin{thm}
\label{main}
Let $u \in C^2 (B_1)$ be a plurisubharmonic function that solves the equation \eqref{CMA}. Suppose
\begin{equation}
\label{c3}
\inf_{B_1} f^{1/n} \geq \lambda > 0, \quad  \text{ and } \quad  \sup_{B_1 } \Delta u  \leq \Lambda.
\end{equation}
Then, there exists a concave function $\wt{F}$ on the space of $2n\times 2n$ real symmetric matrices such that

i) $\wt{F}$ is $\theta$-uniform elliptic, i.e., 
\[
\theta \norm{P} \leq \wt{F} (M +P) - \wt{F} (M) \leq \theta^{-1} \norm{P} , \quad \forall \; M \in \Sym (2n), P\geq 0
\]
where $\theta $ only depends on $\lambda, \Lambda, n$.

ii) $u$ satisfies the equation
\[
\wt{F} (D^2u ) =  f^{1/n} \quad  \text{ in } B_1.
\]   
\end{thm}

Priori to this note, the best result in this direction is obtained by S. Dinew, X. Zhang and X.-W Zhang \cite{Zhang}. They have proved $C^{2,\alpha}$-regularity under the assumption of the $L^{\infty}$-bound of $D^2 u$. Their proof is based on the perturbation argument developed by Trudinger and Wang \cite{Wang}, \cite{Wang1}.

Contrary to the method employed in \cite{Zhang}, we reduce the problem to the uniform elliptic case by constructing an suitable extension of determinant outside a certain set in the space of matrices. The idea of this note is suggested by Prof. Ovidiu Savin.

\bigskip 

\section{The Proof of Thm.\ref{main}}
Let $\Sym(2n)$ be the space of $2n \times 2n$ real symmetric matrices and $\Herm(n)$ be the space of $n \times n$ complex Hermitian matrices. 

Fix the following canonical complex structure 
\[
J = \begin{pmatrix}
0 & -I \\
I & 0
\end{pmatrix} \quad I \text{ is the } n \times n \text{ identity matrix}
\] 
on $\R^{2n}$. Then $\Herm(n)$ can be identified to the subspace
\[
\{M : [M, J] = MJ - JM = 0\} \subset \Sym(2n)
\]
by the map
\[
\imath :  H = A + i B  \mapsto   \begin{pmatrix}
A & - B \\
B & A 
\end{pmatrix}.
\]
In the rest of this note, we always view $\Herm(n)$ as a subspace of $\Sym(2n)$ according to the above identification. 

The complex structure $J$ gives rise the canonical projection $\proj : \Sym(2n) \rightarrow \Herm(n) $
\[
\proj: M \mapsto  \frac{M + J^{t} M J}{2}.
\]

The complex determinant $\det_{\C}$ on Hermitian matrices is related to the real determinant $\det_{\R}$ by
\[
\mathrm{det}_{\R}^{1/2n}[ \proj (M) ] =  \mathrm{det}_{\C}^{1/n}[H] \quad \text{ if } \imath (H) = \proj (M), \; M \in \Sym(2n), \; H \in \Herm(n).
\]

\medskip

We denote
\begin{equation}
\label{notation}
F (M) := \mathrm{det}_{\R}^{1/2n}[ \proj (M) ] ,  \quad   M \in \Sym(2n).
\end{equation}
By the Minkowski inequality, $F$ is a concave function on the set 
\[
\{ M \in \Sym (2n) : \proj (M) > 0 \}.
\]

Now, we give the construction of $\wt{F}$: 
\begin{defn}
\label{MC}
Given $\theta >0$, let $\mathcal{E}_{\theta} \subset \Sym(2n)$ consist of matrices $N$ such that
\[
\theta I  \leq \proj (N)  \leq \theta^{-1} I.
\]
Define, for all $M \in \Sym (2n)$,
\[\begin{split}
\widetilde{F} (M) := \inf\{  \tr [\proj(N) M ]  + c:    i) \tr[\proj( N) X ]+c \geq F(X) \; \forall X \in \mathcal{E}_{\theta} \quad  ii) N  \in \mathcal{E}_{\theta}, c \in \R \}.
\end{split}\]
\end{defn}

\begin{rem}
$\wt{F}$ is the concave envelop of $F$ over the set $\mathcal{E}_{\theta}$.
\end{rem}

\begin{rem}
The above construction is suggessted by Prof. Ovidiu Savin. The author's original approach is to extend level set of $F$ outside $\El_{\theta}$. Though it will give essentially same function as above, the construction in Def.\ref{MC} is more direct and transparent.
\end{rem}

The following lemma is the main ingredient in proving Thm.\ref{main};
\begin{lem}
\label{key}
$\widetilde{F}$ is concave and uniformly elliptic in $\Sym(2n)$, i.e., there exists $\tilde{\theta}>0$ only depends on $\theta$ such that 
\begin{equation}
\label{ue}
\tilde{\theta} \norm{P} \leq  F(M + P) - F(M) \leq \tilde{\theta}^{-1} \norm{P}, \quad \forall M  \in \Sym (2n) , P \geq 0
\end{equation}
Moreover, $\tilde{F}(M) = F (M)$ for all $M \in \mathcal{E}_{\theta}$.
\end{lem}

\begin{proof}
Concavity and agreement on $\mathcal{E}_{\theta}$ follow directly from the construction. We only need to check the ellipticity. Given $M \in \Sym(2n), P \geq 0$, by the definition of $\widetilde{F}$, there are $N_1, N_2 \in \mathcal{E}_{\theta}$ such that
\[
\widetilde{F} (M + P ) = \tr[\proj(N_1) ( M + P)] +c_1
\]
and
\[
\widetilde{F} (M) = \tr[\proj (N_2) (M)] +c_2.
\]
By the minimality, we have
\[
\tr [\proj (N_1) M] +c_1 \geq \tr[\proj (N_2) M] +c_2
\]
and
\[
\tr[\proj(N_1) (M+P) ] +c_1 \leq \tr[\proj (N_2) (M +P)] +c_2.
\]
Then, combine above inequalities, we have
\[
\begin{split}
\widetilde{F} (M + P )  -\widetilde{F} (M) & \geq \tr[\proj(N_1) (M +P)]  - \tr[\proj (N_1) (M)]  \geq \frac{\theta}{4n} \norm{P}
\end{split}
\]
and 
\[
\widetilde{F} (M + P )  -\widetilde{F} (M)  \leq \tr[\proj(N_2) (M +P)] - \tr[\proj (N_2) M] \leq 2n \theta^{-1} \norm{P}.
\]
This completes the proof of the lemma.
\end{proof}

Now, we are ready to complete the proof of Thm.\ref{main}. Since $u \in C^2 (B_1)$ satisfies \eqref{c3}, we have
\[
 \lambda I \leq  \proj (D^2 u ) (x) \leq \Lambda I, \quad \forall  \; x\in B_{1}.
\]
In turn, by taking $\theta = \min\{\lambda, \Lambda^{-1}\}$, we have
\begin{equation}
\label{belong}
D^2 u (x) \in \El_{\theta}, \quad  \forall \; x \in B_{1}.
\end{equation}

Now consider $\wt{F}$ given by the Def.\ref{MC} with respect to $\El_{\theta}$. By \eqref{belong}, Lem.\ref{key} and \eqref{notation}
\[
\wt{F} (D^2u  (x))  =F(D^2 u  (x ))  = f^{1/n} (x) \quad x \in B_{1/2}.
\]
Uniform ellipticity and concavity of $\wt{F}$ have been given in Lem.\ref{key}. This completes the proof of Thm.\ref{main}.

\bigskip 

\textbf{Acknowledgment} The author would like to express his gratitude to Prof. Ovidiu Savin for suggesting the key idea and to Prof. Duong Hong Phong for many inspirational  remarks. The author also would like to thank Valentino Tosatti  and Xiangwen Zhang for many helpful discussion.


\end{document}